\documentclass[10pt]{amsart}
\usepackage{latexsym, amsmath,amssymb}

\usepackage{etoolbox}
\patchcmd{\section}{\scshape}{\bfseries}{}{}
\makeatletter
\renewcommand{\@secnumfont}{\bfseries}
\makeatother

 \numberwithin{equation}{section}

\newtheorem{theorem}{Theorem}

\newtheorem{lemma}[theorem]{Lemma}
\newtheorem{corr}[theorem]{Corollary}

\newtheorem{proposition}[theorem]{Proposition}
\newtheorem{deff}[theorem]{Definition}

\newcommand{\bth}{\begin{theorem}}
\newcommand{\ble}{\begin{lemma}}
\newcommand{\bcor}{\begin{corr}}

\newcommand{\bdeff}{\begin{deff}}

\newcommand{\bprop}{\begin{proposition}}
\newcommand{\ele}{\end{lemma}}
\newcommand{\ecor}{\end{corr}}
\newcommand{\edeff}{\end{deff}}

\newcommand{\eprop}{\end{proposition}}

\newcommand{\la}{\lambda}

\newcommand{\e}{\varepsilon}

\renewcommand{\Pi}{\varPi}

\renewcommand{\epsilon}{\varepsilon}

\begin{document}

\title
{A note on  $L^p$-norms of quasi-modes}

\author{Christopher D. Sogge}
\address{Department of Mathematics,  Johns Hopkins University,
Baltimore, MD 21218}
\email{sogge@jhu.edu}
\author{Steve Zelditch}
\address{Department of Mathematics, Northwestern University, Evanston, IL 60208}
\email{s-zelditch@northwestern.edu}

\thanks{The research of the first author was partially supported by NSF grant \# DMS-1069175 and the second by  \#  DMS-1206527.}

\maketitle
\centerline{\em Dedicated to  Shanzhen Lu on the occasion of his seventy fifth birthday.}
\begin{abstract}
In this note we show how improved $L^p$-estimates for certain types of quasi-modes are naturally equaivalent
to improved operator norms of spectral projection operators associated to shrinking spectral intervals of the appropriate scale.  Using this, one can see that recent estimates that were stated for eigenfunctions also hold for the appropriate types of quasi-modes.
\end{abstract}

\section{Introduction and main results}

Let $(M,g)$ be a compact boundaryless Riemannian manifold of dimension
$n\ge 2$. Let  $0=\la_0<\la_1\le \la_2\le \cdots$ be the frequencies associated with
$\Delta$ and let $e_j(x)$ be an orthonormal basis of  eigenfunctions,
\begin{equation}\label{2}
(\Delta^2+\la_j^2)e_j(x)=0, \quad \text{and } \quad
\int_M|e_j|^2\, dV=1.
\end{equation}
Here $dV=dV_g$ and $\Delta=\Delta_g$ denote the volume element and Laplace-Beltrami operator associated with $(M,g)$, respectively.
Also define the `unit-band'  spectral projection operators (i.e. for frequency intervals of unit length)  by
\begin{equation}\label{3}\chi_{[\la, \la+1]} f(x)=\sum_{\la_j\in [\la,\la+1]}E_jf,\end{equation}
where $E_j$ denotes the projection onto the $j$-th eigenspace, i.e.,
$$E_jf(x)=\langle f,e_j\rangle \, e_j(x),$$
with $\langle \, \cdot \, ,\, \cdot \, \rangle$ being the $L^2(dV)$-inner product.

Many recent works are concerned with $L^p$  norms of eigenfunctions and their relations to the geometry
of the geodesics of $(M, g)$. The first author proved universal bounds on $L^p$ norms of eigenfunctions
\cite{Seig}. For large $p$ the bounds are saturated by zonal (rotationally invariant) spherical harmonics on the standard sphere
$S^n$, which peak at the poles, and for low $p$ they are saturated by elliptic Gaussian beams which concentrate
on closed geodesics (see \cite{Sthesis}). In general, a Riemannian manifold fails to have either kind of eigenfunction and that suggests
that the $L^p$ norms of eigenfunctions are strictly smaller on a typical $(M, g)$ than in the case of the sphere.
The theme of our works \cite{SZ}, \cite{STZ}, \cite{SZfoc} is that  for $p = \infty$ (or for `large p'), the $L^p$
bounds can only be saturated if there exist (partial) poles or self-focal points. Less is known about  necessary
geometric conditions for existence of eigenfunctions saturating the low $L^p$ norms. It is plausible that
that existence of a stable elliptic geodesic is necesssary for saturation of low $L^p$ norms. We refer to
 \cite{Toh},\cite{China}, \cite{Stein} for some recent results and to \cite{SHang}, \cite{Z}, \cite{Z2} for surveys of the problem.

Many of the techniques and results concerning $L^p$ norms  of eigenfunctions (`modes') apply equally well to quasi-modes or approximate eigenfunctions. Quasi-modes are linear combinations of eigenfunctions with eigenvalues concentrated in a short interval.
Specifically, if
$\la$ ranges over a
 sequence $\{\mu_k\}$ tending to $+\infty$, then we say that $\{\psi_{\la}\}$ (with $\la=\mu_k$) is a  quasi-mode of order $s$ if
\begin{equation}\label{1s}
\int_M |\psi_{\la}|^2 \, dV=1, \quad \text{and } \quad
\|(\Delta+\la^2)\psi_{\la}\|_{L^2(M)}= O(\lambda^{-s}).
\end{equation}  Note then that, by our convention of using $\la^2$
instead of $\la$ in the second part of \eqref{1},
the quasi-modes $\psi_{\la}$ are associated with frequencies $\la$. Roughly speaking, the higher the $s$, the
shorter the band of frequencies required to construct the quasi-mode.  We would like to give bounds on $L^p$ norms
of quasi-modes of a given order $s$ and relate them to the geometry of geodesics. For instance, we may ask if
there exist quasi-modes of order $s$ which saturate the $L^p$ bounds and if they resemble zonal spherical harmonics
for high $p$ or Gaussian beams for low p. Of course, quasi-modes peaking at a point or concentrating on a closed
geodesic exist much more often than actual eigenfunctions with these properties. But as in the case of modes,
there are necessary geometric conditions on $(M,g)$ for the construction of quasi-modes of order $s$ which saturate
the $L^p$ bounds of \cite{Seig}. It is not hard to see that there are no conditions if one lets $s = -1$, since
the Schwartz kernel $\chi_{[\lambda, \lambda + 1]}(\cdot, x)$ itself saturates the high $L^p$ bounds and
the spectral  projection $\chi_{[\lambda, \lambda + 1]} (\delta_{\gamma})$ of the delta-function on a closed
geodesic saturates the low $L^p$ bounds.  The question is whether one can find necessary geometric conditions
under which quasi-modes with spectrum in just slightly shrinking intervals can saturate the $L^p$ bounds.
If so, Riemannian manifolds $(M, g)$ which fail to satisfy the conditions are said to have ``improved estimates"
or to satisfy ``sub-convexity bounds". 

As this discussion indicates, closely  related to estimates of modes and quasi-modes are mapping norm estimates of spectral projections
for intervals. If 
\begin{equation}\label{4}
\sigma(p,n)=\begin{cases}
n(\tfrac12-\tfrac1p)-\tfrac12, \quad \tfrac{2(n+1)}{n-1}\le p\le \infty,
\\
\tfrac{n-1}2(\tfrac12-\tfrac1p), \quad  2\le p\le \tfrac{2(n+1)}{n-1},
\end{cases}
\end{equation}
then it was shown in \cite{Seig} that
\begin{equation}\label{5}
\|\chi_{[\la,\la+1]}f\|_{L^p(M)}\le C\la^{\sigma(p,n)}\|f\|_{L^2(M)},
\quad \la \ge 1, \, \, p\ge 2.
\end{equation}
These unit-band spectral projection estimates turn out to be sharp and cannot be improved on {\em any}
compact boundaryless Riemannian manifold, as was shown in \cite{Sbook}.

In this note we consider spectral projections for shrinking intevals. As with quasi-modes, the  question is to find geometric
conditions on $(M, g)$ under which one can improve  or `break convexity' for the mapping norm estimate. 
Breaking convexity bounds on mapping norms of spectral projections means finding a rate  $\epsilon(\lambda)$
of shrinking intervals $[\la,\la+\e(\la)]$ and a condition on $(M, g)$  for which one can prove 
\begin{equation}\label{3.1}
\|\chi_{[\la,\la+\e(\la)]}f\|_{L^p(M)}\le C\sqrt{\e(\la)}
\la^{\sigma(p,n)}\|f\|_{L^2(M)},
\end{equation}
assuming that $\e(\la)>0$ and $\e(\la)$ tends monotonically to zero
as $\la\to \infty$.  Here, as before, $\sigma(p,n)$ is as in \eqref{4}.
The dependence of \eqref{3.1} on $\e(\la)$ is sharp in view of the fact
that \eqref{5} cannot be improved.  To see this one uses orthogonality and
the dual version of \eqref{5}.

 In this
note, we do not give new geometric conditions. Rather our purpose is to show
that the  conditions for improving  mapping norm estimates  or sub-convexity bounds for spectral projection operators associated to shrinking intervals
are equivalent to the conditions  improved bounds for all  quasi-modes of the related order.  The
arguments  are implicit in
\cite{BSSY} and \cite{STZ} but have not previously been made explicit. 
We shall focus on two different types of quasi-modes that naturally occur, although the methods will handle other types as well.

\subsection{Quasi-modes of order $o(\lambda)$}

The first type we call quasi-modes of order $o(\la)$ or
``fat quasi-modes''.  In place of \eqref{1} they satisfy
\begin{equation}\label{1}
\int_M |\psi_{\la}|^2 \, dV=1, \quad \text{and } \quad
\|(\Delta+\la^2)\psi_{\la}\|_{L^2(M)}= o(\lambda).
\end{equation}   
These occur at least implicitly in several recent
works, including \cite{BS}, \cite{bourgain}, \cite{Br},\cite{BSSY},\cite{ChenS}, \cite{Toh}, \cite{STZ}, \cite{SZ} and \cite{SZfoc}.  
We would like to give  conditions on $(M, g)$ which break the convexity bounds, i.e.  prohibit fat quasi-modes from saturating the $L^p$ bounds
of \cite{Seig}. 
Our first result states that breaking convexity bounds for fat quasi-modes is equivalent  to breaking convexity
bounds for  mapping norms of spectral projections for
shrinking intervals. That is, we 
replace the unit intervals $[\la,\la+1]$ by shrinking intervals $[\la,\la+o(1)]$ and  seek sufficient conditions on $(M, g)$ so that
 if $\sigma(p,n)$ is as in \eqref{4}, then
\begin{equation}\label{6}
\|\chi_{[\la,\la+o(1)]}\|_{L^2(M)\to L^p(M)}=o(\la^{\sigma(p,n)}).
\end{equation}
By this we of course mean that for every $\e>0$ there is a $0<\delta <1$
and a $\Lambda<\infty$ (both depending on $\e$) so that
\begin{equation}\label{7}
\|\chi_{[\la,\la+\delta]}\|_{L^2(M)\to L^p(M)}\le \e \la^{\sigma(p,n)},
\quad \la >\Lambda,
\end{equation}
where
$$\chi_{[\la,\la+\delta]}f=\sum_{\la_j\in [\la,\la+\delta]}E_jf,$$
and
$$\|\chi_{[\la,\la+\delta]}\|_{L^2(M)\to L^p(M)}$$
denotes the $L^2(M)\to L^p(M)$ operator norm of $\chi_{[\la,\la+\delta]}$.
Improved bounds \eqref{6}    were shown to be generically true for large exponents
$p$ in \cite{SZ} and \cite{STZ}, and in many cases in \cite{BS} and \cite{Stein} for $2<p<\tfrac{2(n+1)}{n-1}$,

The condition \eqref{1} says that in an $L^2$-sense the spectrum of the
$\psi_{\la}$ are concentrated in intervals of length $o(1)$ about $\la=\mu_k$ as $\mu_k\to \infty$.
In view of this, the following result naturally links $o(\la^{\sigma(p,n)})$ estimates for these $o(\la)$ quasi-modes with the improved $L^p$-bounds in
\eqref{7}.

\begin{theorem}\label{theorem1}
Suppose that $n\le 3$ and $p>2$, $n=4$ and $2<p<\infty$ or $n\ge 5$
and $2<p\le \tfrac{2n}{n-4}$.  Then \eqref{6} is valid if and only if whenever
$\la$ ranges over a sequence $\mu_k\to \infty$ we have
\begin{equation}\label{8}
\|\psi_\la\|_{L^p(M)}=o(\la^{\sigma(p,n)}),
\end{equation}
assuming that the $\psi_\la$ satisfy the $o(\la)$ quasi-mode condition 
\eqref{1}.  If $n=4$ and $p=\infty$ or $n\ge 5$ and $\tfrac{2n}{n-4}<p\le \infty$, the same conclusion is valid if we replace \eqref{1} by the condition
that
\begin{multline}\label{9}
\int_M |\psi_{\la}|^2 \, dV=1, 
\\ \text{and } \quad
\|(\Delta+\la^2)\psi_{\la}\|_{L^2(M)}+\la^{2-n(\tfrac12-\tfrac1p)}
\|\chi_{[2\la,\infty)}\psi_\la\|_{L^p(M)}=o(\la),
\end{multline}
where
$$\chi_{[2\la,\infty)}f=\sum_{\la_j\ge 2\la}E_jf$$
denotes the projection of $f$ onto frequencies $\la_j\in [2\la,\infty)$.
\end{theorem}

Note that for all exponents $p>2$ if $n\le 3$ or for $2<p<\infty$ when $n=4$
or $2<p\le \tfrac{2n}{n-4}$ the second term in the last condition in \eqref{9}
is dominated by the first term by Sobolev embeddings.  Thus, for these
exponents, the condition \eqref{9} is equivalent to the more succinct version
\eqref{1}.  On the other hand for the other exceptional exponents $p>2$
if $n\ge 4$, this term is needed for the conclusion as was argued in 
\cite{STZ}.  This is just due to Sobolev considerations.

\subsection{Quasi-modes of order zero}

The second type of  quasi-mode we consider are those  of order zero, which satisfy
\begin{equation}\label{10}
\int_M|\psi_\la|^2 \, dV=1, \quad 
\text{and } \quad \|(\Delta+\la^2)\psi_\la\|_{L^2(M)}=O(1),
\end{equation}
assuming, as above, that $\la$ ranges over a sequence $\mu_k\to \infty$. We refer to \cite{zworski} for
backgound and references.

Note that the second condition in \eqref{10} says that, in an $L^2$-sense,
the spectrum of the $\psi_\la$ are concentrated in intervals of length
$1/\la$ about $\la$.  Thus, the following result is the natural analog
of Theorem~\ref{theorem1}.

\begin{theorem}\label{theorem2}
Suppose that $n\le 3$ and $p>2$, $n=4$ and $2<p<\infty$ or $n\ge 5$
and $2<p\le \tfrac{2n}{n-4}$.  Then
\begin{equation}\label{11}
\|\chi_{[\la, \la+1/\la]}\|_{L^2(M)\to L^p(M)}=o(\la^{\sigma(p,n)}),
\end{equation}
where $\sigma(p,n)$ is as in \eqref{4} if and only if
whenever
$\la$ ranges over a sequence $\mu_k\to \infty$ we have
\begin{equation}\label{12}
\|\psi_\la\|_{L^p(M)}=o(\la^{\sigma(p,n)}),
\end{equation}
assuming that the $\psi_\la$ satisfy the zero order quasi-mode condition 
\eqref{11}.  If $n=4$ and $p=\infty$ or $n\ge 5$ and $\tfrac{2n}{n-4}<p\le \infty$, the same conclusion is valid if we replace \eqref{11} by the condition
that
\begin{multline}\label{13}
\int_M |\psi_{\la}|^2 \, dV=1, 
\\
 \text{and } \quad
\|(\Delta+\la^2)\psi_{\la}\|_{L^2(M)}+\la^{2-n(\tfrac12-\tfrac1p)}
\|\chi_{[2\la,\infty)}\psi_\la\|_{L^p(M)}=O(1),
\end{multline}
where
$\chi_{[2\la,\infty)}$ is as above.
\end{theorem}

\subsection{A unifying principle}

As we shall see both of these Theorems are a consequence  of the following result, which
is a variation of the ones in \cite{BSSY}.

\begin{proposition}\label{prop}
Suppose that $n\le 3$ and $p>2$, $n=4$ and $2<p<\infty$ or $n\ge 5$
and $2<p\le \tfrac{2n}{n-4}$.  Then
there is a uniform constant $C=C(M,g)$ depending only on our $n$-dimensional Riemannian
manifold $(M,g)$ so that for 
$0<\delta<1$ and sufficiently large $\la$ we have
\begin{multline}\label{14}
\|f\|_{L^p(M)}
\\
\le C\sup_{\mu\in [\la/2,2\la]}\|\chi_{[\mu,\mu+\delta]}\|_{L^2(M)\to
L^p(M)} \times  \bigl(\|f\|_{L^2(M)}+(\delta\la)^{-1}\|(\Delta+\la^2)f\|_{L^2(M)}
\bigr).
\end{multline}
Similarly,
\begin{multline}\label{15}
\|f\|_{L^p(M)}-\|\chi_{[2\la,\infty)}f\|_{L^p(M)}
\\
\le C\sup_{\mu\in [\la/2,2\la]}\|\chi_{[\mu,\mu+\delta]}\|_{L^2(M)\to
L^p(M)} \times  \bigl(\|f\|_{L^2(M)}+(\delta\la)^{-1}\|(\Delta+\la^2)f\|_{L^2(M)}
\bigr),
\end{multline}
for the remaining exponents $p>2$, i.e., $p=\infty$ if $n=4$ and $\tfrac{2n}{n-4}
<p\le \infty$ if $n\ge 5$.
\end{proposition}

\section{Proof that Proposition~\ref{prop} implies  Theorems~\ref{theorem1} and \ref{theorem2}}

We shall start by going over the simple argument which shows how 
Proposition~\ref{prop} implies the two theorems.

Let us first suppose that \eqref{6} is valid.  Thus, given $\e>0$
we can find a fixed $\delta=\delta(\e)\in (0,1)$ so that 
\eqref{7} is valid.  In other words, we have that
$$\sup_{\mu\in [\la/2,2\la]}\|\chi_{[\mu, \mu+\delta]}\|_{L^2(M)\to L^p(M)}
\le \e \la^{\sigma(p,n)},$$
if $\la \gg 1$ is large enough.  By \eqref{14}, we then have that
\begin{equation}\label{2.1}
\|f\|_{L^p(M)}\le C\e \la^{\sigma(p,n)}
\bigl(\|f\|_{L^2(M)}+(\delta\la)^{-1}\|(\Delta+\la^2)f\|_{L^2(M)}\bigr).
\end{equation}
Since we have fixed $\delta\in (0,1)$, if we take $f=\psi_\la$, where
the $\psi_\la$ satisfy \eqref{1}, we conclude that \eqref{2.1} yields
\begin{equation}\label{2.2}
\|\psi_\la\|_{L^p(M)}\le 2C\e \la^{\sigma(p,n)},
\end{equation}
for large enough $\la$.  Thus, since $\e>0$ is arbitrary, we obtain
\eqref{8} assuming that \eqref{6} is valid and $n\le3$ or $2<p<\infty$ if
$n=4$ or $2<p\le \tfrac{2n}{n-4}$ if $n\ge 5$.  Since
\begin{equation}\label{2prime}
n(\tfrac12-\tfrac1p)-1<\sigma(p,n)\end{equation}
for all exponents $p>2$, we similarly obtain \eqref{8}
for the remaining exponents from \eqref{15} assuming that \eqref{6} is valid and that 
the $\psi_\la$ satisfy \eqref{9}.

To prove the remaining half of Theorem~\ref{theorem1}, let
\begin{equation}\label{2.3}
V_{[\la,\la+\delta]}=\bigl\{f\in L^2(M): \, E_jf=0, \, 
\la_j\notin [\la,\la+\delta]\bigr\}
\end{equation}
denote the span of the eigenfunctions whose frequencies lie in
$[\la,\la+\delta]$.  It then follows that
\begin{equation}\label{2.4}
\|\chi_{[\la,\la+\delta]}\|_{L^2(M)\to L^p(M)}=\sup_{f\in V_{[\la,
\la+\delta]}, \, \|f\|_2=1} \, \bigl\|f\bigr\|_{L^p(M)},
\end{equation}
and
\begin{equation}\label{2.5}
\|(\Delta+\la^2)f\|_{L^2(M)}\le 3\delta\la, \quad
\text{if } \, \, f\in V_{[\la, \la+\delta]} \, \, 
\text{and } \, \, \|f\|_{L^2(M)}=1,
\end{equation}
assuming that $0<\delta<1$ and $\la\ge 1$.   Note also that for such
$\delta$, $\lambda$, we have $\chi_{[2\la,\infty)}f=0$ if 
$f\in V_{[\la,\la+\delta]}$.  Thus if \eqref{8} is valid assuming
either \eqref{1} or \eqref{9} as $\la$ ranges over any sequence 
$\mu_k\to \infty$, we conclude that \eqref{6} must hold, completing
the proof of Theorem~\ref{theorem1}.
\medskip

The proof of Theorem~\ref{theorem2} is very similar.  First, if \eqref{11} is
valid then by \eqref{2.1} with $\delta=1/\la$, we conclude using \eqref{14}
again that we must have \eqref{2.2} if the $\psi_\la$ satisfy \eqref{10}
and $n\le 3$ and $p>2$, $2<p<\infty$ if $n=4$ or $2<p\le \tfrac{2n}{n-4}$
if $n\ge 5$.  We similarly obtain \eqref{12} assuming \eqref{13} from
\eqref{15} for the remaining exponents $p>2$, since for these we have
$$n(\tfrac12-\tfrac1p)-2<\sigma(p,n).$$
Since the converse part of Theorem~\ref{theorem2} follows from the corresponding proof of this part of Theorem~\ref{theorem1} (using
$\delta=1/\la$ in \eqref{2.4} and \eqref{2.5}), we obtain
the second theorem as well.


\section{Proof of Proposition~\ref{prop}}

We now prove  Proposition~\ref{prop}.  We first
note that by the Sobolev estimate
$||u||_p \leq  C ||u||_{W^{k, 2}}
$ (with $ k = \frac{n}{2} (\frac{1}{2} - \frac{1}{p}) $)  and orthogonality,  if $2<p<\infty$ and
$n\le 4$ or $2<p\le \tfrac{2n}{n-4}$ and $n\ge 5$, we have
\begin{align}\label{2.6}
\|\chi_{[2\la,\infty)}f\|_{L^p(M)}&\le C\|(-\Delta)^{\frac{n}2(\frac12-\frac1p)} \chi_{[2\la,\infty)}f\|_{L^2(M)}
\\
&\le C\la^{n(\frac12-\frac1p)-2}\|(\Delta+\la^2)\chi_{[2\la,\infty)}f\|_{L^2(M)} \notag
\\
&\le C\la^{n(\frac12-\frac1p)-2}\|(\Delta+\la^2)f\|_{L^2(M)}.\notag
\end{align}
Similarly, for every $\e>0$ we have
\begin{equation}\label{2.7}
\|\chi_{[2\la,\infty)}f\|_{L^\infty(M)}
\le C_\e \la^{\frac{n}2-2+\e}\|(\Delta+\la^2)f\|_{L^2(M)},
\quad \text{if } \, \, n\le 3.
\end{equation}
Also, as we noted before, \eqref{5} cannot be improved, and, therefore, there must be a uniform constant $c>0$, 
which is independent of $0<\delta<1$, so that if $\la$ is large enough we have
$$c\delta \la^{\sigma(p,n)}\le \sup_{\mu\in [\la/2,2\la]}\|
\chi_{[\mu,\mu+\delta]}\|_{L^2(M)\to L^p(M)},$$
and so
\begin{multline}\label{2.8}
\la^{\sigma(p,n)-1}\|(\Delta+\la^2)f\|_{L^2(M)}
\\
\le C\sup_{\mu\in [\la/2,2\la]}\|\chi_{[\mu,\mu+\delta]}\|_{L^2(M)\to
L^p(M)} \times (\delta\la)^{-1}\|(\Delta+\la^2)f\|_{L^2(M)}.
\end{multline}
Therefore, if we recall \eqref{2prime} and combine \eqref{2.6}-\eqref{2.8}, we conclude that
if $2<p\le \infty$ and $n\le 3$ or $2<p<\infty$ and $n=3$ or
$2<p\le \tfrac{2n}{n-4}$ and $n\ge 5$, then there is a uniform
constant $C$, which is independent of $0<\delta<1$, so that for 
sufficiently large $\la$ we have
\begin{equation*}
\|\chi_{[2\la,\infty)}f\|_{L^p(M)}
\le C(\delta\la)^{-1}\|(\Delta+\la^2)f\|_{L^2(M)}
\times \sup_{\mu\in [\la/2,2\la]}
\|\chi_{[\mu,\mu+\delta]}\|_{L^2(M)\to L^p(M)}.
\end{equation*}

From this, we deduce that we would have \eqref{14} and \eqref{15} if we 
could show that 
\begin{multline}\label{2.9}
\|\chi_{[0,2\la)}f\|_{L^p(M)}
\\
\le C\sup_{\mu\in [\la/2,2\la]}\|\chi_{[\mu,\mu+\delta]}\|_{L^2(M)\to
L^p(M)} \times  \bigl(\|f\|_{L^2(M)}+(\delta\la)^{-1}\|(\Delta+\la^2)f\|_{L^2(M)}
\bigr),
\end{multline}
for all exponents $2<p\le \infty$ if $n\ge 2$ and $\delta$ and $\la$ are as above.

To prove this we write
\begin{equation}\label{2.11}
\chi_{[0,2\la)}f=\chi_{[\la-\delta,\la+\delta]}f
+\sum_{k=1}^\infty \, \sum_{|\la-\la_j|\in (\delta k,\delta(k+1)]}
\chi_{(0,2\la)}E_jf.
\end{equation}
Clearly
\begin{align}\label{2.12}
\|\chi_{[\la-\delta,\la+\delta]}f\|_{L^p(M)}
&\le \|\chi_{[\la-\delta,\la+\delta]}\|_{L^2(M)\to L^p(M)}\|f\|_{L^2(M)}
\\
&\le  2\sup_{\mu\in [\la/2,2\la]}\|\chi_{[\mu,\mu+\delta]}\|_{L^2(M)\to
L^p(M)}\times \|f\|_{L^2(M)}. \notag   
\end{align}
Similarly,
\begin{align*}
\bigl\|&\sum_{|\la-\la_j|\in (\delta k,\delta(k+1)]}
\chi_{(0,2\la)}E_jf\bigr\|_{L^p(M)}
\\
&\le \sup_{\mu\in [\la/2,2\la]}\|\chi_{[\mu,\mu+\delta]}\|_{L^2(M)\to
L^p(M)} \times  
\bigl\|\sum_{|\la-\la_j|\in (\delta k,\delta(k+1)]}
\chi_{(0,2\la)}E_jf\bigr\|_{L^2(M)} \notag
\\
&\le (\delta k\la)^{-1}\sup_{\mu\in [\la/2,2\la]}
\|\chi_{[\mu,\mu+\delta]}\|_{L^2(M)\to
L^p(M)} \times  
\bigl\|\sum_{|\la-\la_j|\in (\delta k,\delta(k+1)]}
(\Delta+\la^2)E_jf\bigr\|_{L^2(M)}. \notag  
\end{align*}
Therefore, by orthogonality and an application of the Cauchy-Schwarz
inequality, we deduce that the $L^p(M)$-norm of the the last term
in \eqref{2.11} is
$$\le C\sup_{\mu\in [\la/2,2\la]}\|\chi_{[\mu,\mu+\delta]}\|_{L^2(M)\to
L^p(M)} \times (\delta\la)^{-1}\|(\Delta+\la^2)f\|_{L^2(M)}.
$$
From this and \eqref{2.12}, we obtain \eqref{2.9}, which completes
the proof of Proposition~\ref{prop}.

\section{Applications to breaking convexity bounds}

In this section, we survey some recent results on breaking convexity bounds 
on mapping norms for $\chi_{[\la,\la+\e(\la)]}$  of the form \eqref{3.1} . 
By  Proposition~\ref{prop},  breaking such convexity bounds is  equivalent to proving estimates
\begin{equation}\label{3.2}\|\psi_\la\|_{L^p(M)}\le C\sqrt{\e(\la)}\la^{\sigma(p,n)},
\end{equation}
 for quasi-modes satisfying
\begin{equation}\label{3.3}
\|\psi_\la\|_{L^2(M)}\le 1
\quad \text{and }\quad \|(\Delta+\la^2)\psi_\la\|_{L^2(M)}\le \e(\la) \la,
\end{equation}
assuming that $2<p\le \infty$ if $n\le 3$, $2<p<\infty$ if $n=4$ and
$2<p\le \tfrac{2n}{n-3}$ if $n\ge 5$.  For the remaining exponents $p>2$
if $n=4$ or $n\ge 5$, we can obtain \eqref{3.2} if we also assume that
\begin{equation}\label{3.4}
\la^{2-n(\frac12-\frac1p)}\|\chi_{[2\la,\infty)}\psi_\la\|_{L^p(M)}=O(1).
\end{equation}
Note that the second condition in \eqref{3.3} is natural for \eqref{3.1} since
it says that, in an $L^2$-sense, most of the spectrum of the $\psi_\la$ 
are concentrated in intervals of size $O(\e(\la)\la)$ about $\la$.

\subsection{Nondissipative focal points exist if  convexity cannot be broken}

It is proved in \cite[Lemma 6.1]{SZ} and \cite[Theorem 2]{STZ} that  for generic $(M,g)$ there exists a monotonically
decreasing function $\e(\la)\to 0$ (the rate depending on $(M,g)$) so that 
\eqref{3.1} is valid.  Thus, in this case, we also get \eqref{3.2} for
$p=\infty$ and the above type of quasi-modes. However, $\epsilon(\lambda)$ is not described explicitly.

If one {\it cannot} break the convexity bound on mapping norms, then one has
\begin{equation}\label{3.5}
\|\chi_{[\la+o(1)]}\|_{L^2(M)\to L^\infty(M)}=\Omega(\la^{\frac{n-1}2}),
\end{equation}
by which we mean that there is a uniform constant $c>0$ so that
\begin{equation}\label{3.6}
\limsup_{\la\to \infty}\la^{-\frac{n-1}2}\|\chi_{[\la,\la+\delta]}\|_{L^2(M)
\to L^\infty(M)}\ge c,
\end{equation}
whenever $\delta>0$.  In the real analytic case in \cite{SZfoc} we found a 
necessary and sufficient condition for this.  Specifically, we showed that 
\eqref{3.5} is valid if and only if our real analytic manifold $(M,g)$
has a nondissipative self-focal point $p$. `Self-focal' means that there exists a time $T(p)$ so that,  for all $\xi \in S^*_p M$ (the unit tangent space),
the geodesic $\gamma_{p, \xi}(t)$ with initial data $(p, \xi)$ returns to $p$ at time $T(p)$.  The geodesic is not
assumed to smoothly close up (i.e. to be a closed geodesic), but just to loop back. The map $\Phi_p: S^*_p M \to S^*_p M$
taking the initial direction $\xi$ to the terminal direction $\gamma_{x, \xi}'(T(p))$ is known as the first return map.
It induces a unitary operator $U_p: L^2(S^*_p M) \to L^2(S^*_p M)$ defined by $$U_p(f) = f(\Phi_p(\xi)) \sqrt{J_p(\xi)}, $$
where $\Phi_p^* |d\xi| = J_p(\xi) |d \xi|$ is the Jacobian with respect to the induced Riemannian volume density on $S^*_p M$. 
We say that a self-focal point is dissipative if $U_p$ has no $L^2$ eigenfunction. Otherwise it is nondissipative.

To obtain \eqref{3.5} we showed that if there exists a nondissipative self-focal point $x_0$, then
if $0<\delta<1$ we have
\begin{equation}\label{3.7}
\limsup_{\la\to \infty}\la^{-(n-1)}\sum_{j=0}^\infty \rho(\delta^{-1}(\la
-\la_j)) |e_j(x_0)|^2\ge c_0,
\end{equation}
for some $c_0>0$. Here,  $\rho\in {\mathcal S}({\mathbb R})$ is any test function satisfying
$$\rho, \Hat \rho\ge 0, \, \, \,
\text{supp }\Hat \rho \subset [-1,1], \, \, \,
\rho(s)\ge 1, \, \, \text{if } \, |s|\le 1.$$  
We prove \eqref{3.7} by using the von Neumann ergodic theorem for $U_p$ together
with  arguments adapted from  \cite{DG}, \cite{Sa}
and \cite{SV}.
We  obtain \eqref{3.6} from
\eqref{3.7} due to the fact that
$$\|\chi_{[\la,\la+\delta]}\|^2_{L^2(M)\to L^\infty(M)}
=\sup_{x\in M}\sum_{\la_j\in [\la,\la+\delta]}|e_j(x)|^2.
$$
Conversely, if \eqref{3.6} is valid one of course gets that, for every
$0<\delta<1$,
$$\limsup_{\la\to \infty}\Bigl( \, \sup_{x\in M}\la^{-(n-1)}
\sum_{j=0}^\infty \rho(\delta^{-1}(\la-\la_j))|e_j(x)|^2\Bigr)\ge c_0>0.
$$

It  follows from \eqref{3.5} and  Proposition \ref{prop}  that  if our real
analytic $(M,g)$ has a nondissipative self-focal point then there is
a $o(\la)$ sequence of quasi-modes satisfying
\begin{equation}\label{3.8}
\|\psi_\la \|_{L^\infty(M)}=\Omega(\la^{\frac{n-1}2}).
\end{equation} We proved this in \cite{SZfoc}  by using a special case of Proposition \ref{prop}. 

\subsection{Existence of  a quasi-mode of order zero with maximal growth of $L^{\infty}$ norms?}

A natural question 
 in \cite{SZfoc} is whether  \eqref{3.5} implies existence of a   sequence of quasi-modes of order zero which
saturate the $L^{\infty}$ bounds. That is, if one has an $L^2$ eigenfunction of $U_p$ at a self-focal point p,
can one construct from it a quasi-mode of order zero with maximal sup-norm growth?  This would definitely
be the case if the  $L^2$ eigenfunction were $C^{\infty}$. In that case, we could use it as the symbol for
a quasi-mode of order zero. In \cite{STZ} such a construction was made 
in the special case of a self-focal point for which $\Phi_p = Id$, i.e. for which all of the loops at $p$ were
smoothly closed. 

  In
view of Theorem~\ref{theorem2},   existence of a zeroth order quasi-mode of maximal sup norm growth would  imply that
\begin{equation}\label{3.9}
\limsup_{\la\to \infty} \la^{-\frac{n-1}2}\|\chi_{[\la,\la+1/\la]}
\|_{L^2(M)\to L^\infty(M)}>0.
\end{equation}


\subsection{Logarithmic shrinking of intervals for non-positively curved manifolds}
In the case of non-positively curved manifolds we can be more specific about the shrinking rate of
intervals and of $L^p$ estimates.

Implicit in \cite{Berard} is the fact that \eqref{3.1} is valid with
$\e(\la)=1/\ln(\la)$ if $(M,g)$ has nonpositive sectional
curvatures (see \cite[\S 3.6]{SHang}).  This result was extended
to all exponents $\tfrac{2(n+1)}{n-1}<p\le \infty$ in 
\cite{HT} (see also \cite[\S 4]{BSSY}).  Thus, for such manifolds we have
\begin{equation} \label{psi} \|\psi_\la\|_{L^p(M)}\le C\la^{\sigma(p,n)}/\sqrt{\ln \la}, \quad  
\frac{2(n+1)}{n-1}<p\le \infty, \end{equation}
assuming that 
$$\|\psi_\la\|_{L^2(M)}\le 1 \quad \text{and } \quad
\|(\Delta+\la^2)\psi_\la\|_{L^2(M)}\le C\la/\ln \la,$$
as well as \eqref{3.4} if $n=4$ or $n\ge 5$ and $p$ is large.

The logarithmic
 estimate \eqref{psi}  reflects the hyperbolicity of the geodesic flow on negatively curved manifolds. It is obtained
more generally in microlocal analysis of quasimodes around hyperbolic closed geodesics. It is plausible that
one does not need global curvature assumptions but just the assumption that all closed geodesics of
$(M, g)$ are hyperbolic.  In the case of {\it quantum integrable} Laplacians (such as a surface of revolution
shaped like a peanut), it is possible to construct eigenfunctions
which concentrate microlocally on a hyperbolic closed geodesic and whose $L^p$ norms in tubes around the closed
geodesic  saturate the inequality \eqref{psi}. 
 We refer to \cite[\S 5.1]{Z2}   for further discussion and references.

 It is  possible to construct
quasi-modes $\psi_{\lambda}$ on negatively curved surfaces which partially concentrate on hyperbolic closed geodesics. A quasi-mode satisfying \eqref{3.3} with 
$\e(\la)=1/\ln(\la)$ has been constructed recently by S. Brooks \cite{Br}  on a compact hyperbolic surface whose microlocal
lift (or defect measure)  has mass $\geq \delta $ for a certain  $\delta > 0$  along a closed geodesic $\gamma$. The
estimate of \eqref{psi} gives an upper bound on its $L^p$ norms, and for $p = 6 + \epsilon$ is precisely the order of magnitude of a quasi-mode associated
to a hyperbolic geodesic. It would be interested to see if Brooks' quasi-mode saturates the bound \eqref{psi}. 

\subsection{$\epsilon(\lambda) = \la^{-\e(p,n)}$ on a flat torus}

In the case of the standard $n$-torus ${\mathbb T}^n={\mathbb R}^n/{\mathbb Z}^n$,
it was also shown in \cite{BSSY}  that we have estimates of the form
$$\|\chi_{[\la, \la +\la^{-\e(p,n)}]}\|_{L^2({\mathbb T}^n)
\to L^p({\mathbb T}^n)}\le \la^{\sigma(p,n)-\e(p,n)/2}, 
\quad \frac{2(n+1)}{n-1}<p\le \infty,$$
for certain powers $\e(p,n)>0$.  Therefore, by the above, we have
$$\|\psi_\la\|_{L^p({\mathbb T}^n)}\le C\la^{\sigma(p,n)-\e(p,n)/2},
\quad \frac{2(n+1)}{n-1}<p\le \infty,$$
assuming that
$$\|\psi_\la\|_{L^2({\mathbb T}^n)}\le 1 \quad
\text{and }\quad \|(\Delta+\la^2)\psi_\la\|_{L^2({\mathbb T}^n)}\le 
C\la^{1-\e(p,n)},$$
as well as \eqref{3.4} if $p$ is large enough and $n\ge 4$.

\end{document}